\documentclass[12pt]{article}
\title{Free Arrangements over Finite Field}
\author{Masahiko Yoshinaga \thanks{International Centre for Theoretical Physics, 
Strada Costiera 11, Trieste 34014, Italy, 
e-mail address: myoshina@ictp.it}}
\date{May 31, 2006}

\usepackage{amssymb}

\newtheorem{Def}{Definition}

\newtheorem{Thm}[Def]{Theorem}
\newtheorem{Lemma}[Def]{Lemma}
\newtheorem{Cor}[Def]{Corollary}
\newtheorem{Rem}[Def]{Remark}

\newtheorem{Example}[Def]{Example}

\newcommand{\bbF}{\mathbb{F}}

\newcommand{\bbK}{\mathbb{K}}

\newcommand{\bbQ}{\mathbb{Q}}

\newcommand{\bbZ}{\mathbb{Z}}

\newcommand{\calA}{\mathcal{A}}

\newcommand{\rmDer}{\mathrm{Der}}

\newcommand{\owari}{\hfill$\square$}

\begin{document}
\maketitle

\begin{abstract}
The freeness of hyperplane arrangements in a three dimensional 
vector space over finite field is 
discussed. We prove that if the number of hyperplanes is 
greater than some bound, then the freeness is determined by 
the characteristic polynomial. 
\end{abstract}

\section{Introduction}
A hyperplane arrangement $\calA$ in an $\ell$-dimensional vector space 
$V$ is a finite set of linear subspaces of codimension one. 
An arrangement 
$\calA$ is said to be free when the associated module of logarithmic 
vector fields is a free module. 
Study of arrangements of this class was started by 
K. Saito \cite{sai-log} 
and 
a remarkable factorization 
theorem was proved by H. Terao \cite{ter-free}. 
This theorem asserts that the characteristic polynomial 
$\chi(\calA, t)$ 
of a free arrangement completely factors into linear 
polynomials in $\bbZ$. 
It imposes a necessary condition on the structure of 
intersection lattice 
$L(\calA)$ for an arrangement $\calA$ to be free. 
The Terao conjecture is the problem 
to ask the converse: does the structure 
of $L(\calA)$ characterize freeness of $\calA$? 
This conjecture is still open even in the case $\ell=3$. 
The purpose of this paper is to propose an affirmative 
result over finite field $\bbF_q$ in the case $\ell=3$. 
Our main result asserts that if the number of hyperplanes satisfies 
$|\calA|\geq 2q-2$, then $\calA$ is free exactly when 
$\chi(\calA, q)\chi(\calA, q-1)=0$. 
Proofs are based on Terao's addition-deletion theorem 
(Theorem \ref{thm:addel}) and Crapo-Rota's method of counting 
points by using characteristic polynomials 
(Theorem \ref{thm:1}, Theorem \ref{thm:cra-rot}).

\section{Freeness and characteristic polynomials}

Let $V$ be an $\ell$-dimensional linear space over a field 
$\bbK$ 
and $S:=\bbK[V]$ be the algebra of polynomial functions on $V$ that is 
naturally isomorphic to $\bbK[z_1, z_2, \cdots, z_\ell]$ for 
any choice of basis $(z_1, \cdots, z_\ell)$ of $V^*$. 
A (central) hyperplane arrangement $\calA$ is a finite 
collection of one-codimensional linear subspaces 
($=$hyperplanes) in $V$. 
For each hyperplane $H$ of $\calA$, 
fix 
a nonzero linear form $\alpha_H\in V^*$ vanishing on $H$ and 
put $Q:=\prod_{H\in\calA}\alpha_H$. 

The characteristic polynomial of $\calA$ is defined as 
$$
\chi(\calA, t)=\sum_{X\in L(\calA)}\mu(X)t^{\dim X},
$$
where $L(\calA)$ is the lattice consisting of intersections of 
elements of $\calA$, ordered by reverse inclusion, 
$\hat{0}:=V$ is the unique minimal element of $L(\calA)$ and 
$\mu:L(\calA)\longrightarrow\bbZ$ is the M\"obius function 
defined as follows: 
\begin{eqnarray*}
\mu(\hat{0})&=&1,\\
\mu(X)&=&-\sum_{Y<X}\mu(Y),\ \mbox{if}\ \hat{0}<X.
\end{eqnarray*}
Fix a hyperplane $H\in\calA$, we obtain two associated arrangements: 
deletion $\calA'=\calA\backslash\{H\}$ and restriction $\calA''=H\cap
\calA'$. 
The characteristic polynomials for these arrangements satisfy 
the following inductive formula 
$$
\chi(\calA, t)=\chi(\calA', t)-\chi(\calA'', t).
$$

Denote by $\rmDer_V$ the $S$-module of all polynomial 
vector fields over $V$. For a given arrangement $\calA$, 
we define the module of logarithmic vector fields as 
$$
D(\calA)=\{ \delta\in\rmDer_V\ |\ \delta(\alpha_H)\in \alpha_H S,\
\forall H\in\calA\}.
$$
An arrangement $\calA$ is said to be {\bf free}, if 
$D(\calA)$ is a free $S$-module, and then the multiset of 
degrees $\exp(\calA):=(d_1, d_2, \cdots, d_\ell)$ 
of homogeneous basis of $D(\calA)$ 
is called the {\bf exponents}. 
The following theorems are due to H. Terao. 

\begin{Thm}
\normalfont\label{thm:addel}
(\cite{ter-arr}, {\cite[Thm 4.52]{orl-ter}})\\
Let $\calA$ be a non-empty arrangement in $\bbK^3$. Let 
$(\calA, \calA', \calA'')$ be a triple as above. Then, 
any two of the following imply the third: 
\begin{itemize}
\item $\calA$ is free with exponents $(d_1, d_2, d_3)$. 
\item $\calA'$ is free with exponents $(d_1, d_2, d_3-1)$. 
\item $|\calA''|=d_1+d_2$. 
\end{itemize}
\end{Thm}

\begin{Thm}
\normalfont
\cite{ter-free}
If $\calA$ is a free arrangement with exponents $(d_1, d_2, \cdots, 
d_\ell)$ 
then the characteristic polynomial is 
$$
\chi(\calA, t)=(t-d_1)\cdots(t-d_\ell).
$$
\end{Thm}

\begin{Example}
\label{ex:all}
Let $V$ be an $\ell$ dimensional vector field over a finite field 
$\bbF_q$ of $q$ elements and $\calA_{all}(V)$ be the collection of 
all hyperplanes in $V$. Put 
$$
\delta_k=\sum_{i=1}^{\ell}x_i^{q^k}\frac{\partial}{\partial x_i}. 
$$
Then for any linear form $\alpha$, $\delta_k\alpha=\alpha^k$. 
Hence $\delta_k\in D(\calA_{all})$. From Saito's criterion \cite{sai-log}, 
$\delta_0, \delta_1, \cdots, \delta_{\ell-1}$ 
form a basis of $D(\calA_{all})$ and 
exponents are $(1, q, \cdots, q^{\ell-1})$. 
\end{Example}


\section{Arrangements over finite fields}

If $\bbK$ is a finite field, the characteristic 
polynomial $\chi(\calA, t)$ has a special meaning. 
The following theorem and its proof is found in \cite[2.69]{orl-ter}, 
and it is a special version of more general result 
obtained in \cite{cra-rot} (see also 
Theorem \ref{thm:cra-rot} below). 

\begin{Thm}
\label{thm:1}
Let $\calA$ be a hyperplane arrangement in 
$V\cong \bbF_q^\ell$. Let $|M(\calA)|$ denote the 
cardinality of the complement. Then 
$$
|M(\calA)|=\chi(\calA, q).
$$
\end{Thm}

This theorem has successfully applied by Athanasiadis \cite{ath-cha} 
to compute characteristic polynomials for arrangements defined 
over $\bbQ$. 

An arrangement $\calA$ over a field $\bbF_q$ can be naturally 
considered as an arrangement over an extended field $\bbF_{q^k}$. Since 
field extensions do not change the intersection lattice, 
the characteristic polynomial is also unchanged. 
Hence from Theorem \ref{thm:1}, we have 
$$
\left|
V\otimes_{\bbF_q}{\bbF_{q^k}}-\bigcup_{H\in\calA}(H\otimes{\bbF_{q^k}})
\right|
=\chi(\calA, q^k).
$$
Not only $\chi(\calA, q)$, but also $\chi(\calA, q^k)$ expresses 
the cardinality of points of complement of 
arrangement $\calA\otimes\bbF_{q^k}$. 
We have the following lemma. 

\begin{Lemma}
\label{lem:easy}
Let $\calA$, $\calA_1$ and $\calA_2$ be arrangements 
in $\bbF_q^\ell$, then 
\begin{itemize}
\item[(i)] $\chi(\calA, q^k)\geq 0$ for all $k\in\bbZ_{>0}$. 
\item[(ii)] If $\calA_1\subset \calA_2$, then 
$\chi(\calA_1, q^k)\geq \chi(\calA_2, q^k)$ for all $k\in\bbZ_{>0}$. 
\end{itemize}
\end{Lemma}

The following Theorem is due to H. Crapo and G. -C. Rota. 
It contains Theorem \ref{thm:1} as a special case ($k=1$). 
Here we deduce from Theorem \ref{thm:1}. 

\begin{Thm}
\cite[\S 16, Theorem 1]{cra-rot}
\label{thm:cra-rot}
Let $\calA$ be an arrangement in $\bbF_q^\ell$. 
The number of ordered points $(p_1, \cdots, p_k)\in(\bbF_q^\ell)^k$  
satisfying the following condition (*) 
is $\chi(\calA, q^k)$: 
\begin{itemize}
\item[(*)] For each hyperplane $H\in\calA$, there exists 
at least one point $p_i$ such that $p_i\notin H$. 
\end{itemize}
\end{Thm}

\noindent
{\bf Proof.}
Recall that $\bbF_{q^k}$ is an $k$-dimensional vector space 
over $\bbF_q$. Let $x_1, \cdots, x_k$ be a $\bbF_q$-basis of $\bbF_{q^k}$. 
Then the point $P$ in $\bbF_{q^k}^\ell$ is expressed as 
$$
P=\left(\sum_{j=1}^{k}a_{1j}x_j, \sum_{j=1}^{k}a_{2j}x_j, \cdots, 
\sum_{j=1}^{k}a_{\ell j}x_j \right), 
$$
where $a_{ij}\in\bbF_q$. Since the defining equation of 
$H\otimes \bbF_{q^k}$ is a linear form with coefficients in $\bbF_q$, 
the point $P$ is contained in $H\otimes \bbF_{q^k}$ if and only if 
$$
(a_{1j}, a_{2j}, \cdots, a_{\ell j})\in H,\ \forall j=1, \cdots, k. 
$$
Hence $P$ is in the complement of arrangement $\calA\otimes\bbF_{q^n}$ 
if and only if for each $H\in\calA$, there exists at least one 
$j\in\{1, \cdots, k\}$ such that 
$(a_{1j}, a_{2j}, \cdots, a_{\ell j})\notin H$. 
So this gives a bijection between complement of $\calA\otimes\bbF_{q^k}$ 
and ordered $k$ points in $\bbF_q^\ell$ satisfying (*). 
\owari

The next result is shown immediately from Theorem \ref{thm:cra-rot}, 
but we give an alternative proof, since the arguments used in the 
proof is prototypical for later. 

\begin{Lemma}
Let $\calA$ be an arrangement in $\bbF_q^\ell$. If the characteristic 
polynomial satisfies 
$\chi(\calA, q^k)=0$ for some $k$, 
then $\chi(\calA, q^j)=0$ for $0\leq j\leq k$. 
\end{Lemma}

\noindent
{\bf Proof.}
The proof is done by induction on the dimension $\ell$ and 
``descending'' induction on the number $|\calA|$ of 
hyperplanes. If $|\calA|$ is maximal, in other words, 
$\calA=\calA_{all}(\bbF_q^\ell)$, then the characteristic 
polynomial is 
$$
\chi(\calA_{all}(\bbF_q^\ell), t)=(t-1)(t-q)\cdots (t-q^{\ell-1}), 
$$
so the lemma is trivial. 
In the case $\ell =2$ is also trivial. 
In general, let $\calA'$ be an arrangement such that 
$\calA'\subsetneq \calA_{all}(\bbF_q^\ell)$ and assume 
$\chi(\calA', q^k)=0$. We can find a 
hyperplane $H$ which is not contained in $\calA'$ and define 
$\calA=\calA'\cup\{H\}$. 
From Lemma \ref{lem:easy}, we have 
$\chi(\calA', q^k)\geq\chi(\calA, q^k)=0$. 
From the inductive hypothesis on the number of hyperplanes, 
we have 
$$
\chi(\calA, q^j)=0, \mbox{ for }0\leq j\leq k. 
$$
Denote $\calA''$ the restriction $H\cap \calA'$, we have 
$\chi(\calA'', q^k)=\chi(\calA', q^k)-\chi(\calA, q^k)=0$. 
Since $\dim H <\ell$, we have 
$$
\chi(\calA'', q^j)=0, \mbox{ for }0\leq j\leq n.
$$
Again from inductive formula, we have 
$\chi(\calA', q^j)=0, \mbox{ for }0\leq j\leq n$. 
\owari

Using the above lemma, we can characterize $\calA_{all}(\bbF_q^\ell)$. 

\begin{Cor}
Let $\calA$ be an arrangement in $\bbF_q^\ell$. The following 
conditions are equivalent. 
\begin{itemize}
\item[(a)] $\calA=\calA_{all}(\bbF_q^\ell)$. 
\item[(b)] $|\calA|=\frac{q^\ell -1}{q-1}$ 
\item[(c)] $\chi(\calA, t)=(t-1)(t-q)\cdots(t-q^{\ell-1})$ 
\item[(d)] $\chi(\calA, q^{\ell-1})=0$. 
\end{itemize}
\end{Cor}

\noindent
{\bf Proof.}
(a)$\Leftrightarrow$(b)$\Leftrightarrow$(c)$\Rightarrow$(d) is trivial. 
(d)$\Rightarrow$(c) is from the above lemma. 
\owari

Here we assume $\ell=3$, and give some combinatorial characterization 
for freeness. 

\begin{Lemma}
\label{lem:ineq}
Let $\calA$ be an essential free arrangement in $\bbF_q^3$ with  
exponents $(1, d_2, d_3)$. If $d_2\leq d_3$, then $d_2\leq q$. 
\end{Lemma}

\noindent
{\bf Proof.}
Note that $d_2$ is the minimal degree of the logarithmic 
vector field $\delta\in D(\calA)$ which does not a 
polynomial multiple of the Euler vector field 
$\delta_0$ (see Example \ref{ex:all} for the notation). 
Since $\delta_1$ is contained in $D(\calA)$, such a minimal 
degree can not be greater than $q=\deg \delta_1$. 
\owari

\begin{Thm}
\label{thm:2q}
Let $\calA$ be an arrangement in $\bbF_q^3$. \\
(1) If $\chi(\calA, q)=0$, then $\calA$ is free with 
exponents $(1, q, |\calA|-q-1)$. \\
(2) If $|\calA|\geq 2q$ and $\calA$ is free then $\chi(\calA, q)=0$.
\end{Thm}

\noindent
{\bf Proof.}
(1) Suppose that $\chi(\calA, q)=0$. We prove the freeness 
by descending induction on $|\calA|$. Recall that $\calA=\calA_{all}$ 
is free with exponents $(1, q, q^2)$ (see Example \ref{ex:all}). 
In general, choose a hyperplane $H$ which is not a member of $\calA$. 
Then by Lemma \ref{lem:easy}, we have 
$$
0=\chi(\calA, q)\geq \chi(\calA\cup\{H\}, q)\geq 0, 
$$
hence $\chi(\calA\cup\{H\}, q)=0$. 
From the inductive hypothesis, 
$\calA\cup\{H\}$ is free with exponents $(1, q, |\calA|-q)$. 
Then Theorem \ref{thm:addel} enable us to conclude that $\calA$ is free. 

(2) Suppose $\calA$ is free. Note that from the assumption on $|\calA|$, 
$\calA$ is an essential arrangement. Hence the characteristic 
polynomial is of the form $\chi(\calA, t)=(t-1)(t-d_2)(t-d_3)$ with 
integers $d_2\leq d_3$ which satisfy $d_2+d_3=|\calA|-1\geq 2q-1$. 
From the Lemma \ref{lem:ineq}, 
we have $d_2\leq q\leq d_3$. However $d_2<q<d_3$ contradicts 
the Lemma \ref{lem:easy} (i) $\chi(\calA, q)\geq 0$, 
we have either 
$$
\begin{array}{ll}
(d_2, d_3)=(q, |\calA|-q-1) & \mbox{ if } |\calA|>2q, \mbox{ or }\\
(d_2, d_3)=(q-1, q) & \mbox{ if } |\calA|=2q.
\end{array}
$$

\owari

By a similar argument to (1), 
we have the following theorem for higher dimensional 
cases. 
\begin{Thm}
Let $\calA$ be an arrangement in $\bbF_q^\ell$. If $\calA$ satisfies
$$
\chi(\calA, q^{\ell-2})=0,
$$
then $\calA$ is free with exponents $(1, q, \cdots, q^{\ell-2},
|\calA|-1-q-\cdots-q^{\ell-2})$.
\end{Thm}

\begin{Rem}
\normalfont
The argument used in the proof of Theorem \ref{thm:2q} (1) can be 
considered as an example of ``supersolvable resolution'' in \cite{zie-mat}. 
\end{Rem}

In the next result we will treat the cases 
$|\calA|=2q-1$ and $2q-2$. 

\begin{Thm}
Suppose that $\chi(\calA, q)\neq 0$. \\
(1) When $|\calA|=2q-1$, $\calA$ is free if and only if 
$\chi(\calA, t)=(t-1)(t-q+1)^2$. \\
(2) When $|\calA|=2q-2$, $\calA$ is free if and only if 
$\chi(\calA, t)=(t-1)(t-q+1)(t-q+2)$. 
\end{Thm}

\noindent
{\bf Proof.}
(1) The similar arguments above from the fact 
$\chi(\calA, q)\geq 0$ shows that 
the freeness of $\calA$ implies $\chi(\calA, t)=(t-1)(t-q+1)^2$. 
Conversely, suppose $\chi(\calA, t)=(t-1)(t-q+1)^2$, then 
$$
|M(\calA)|=\chi(\calA, q)=q-1. 
$$
This means that there exists a line $L\subset \bbF_q^3$ such that 
$$
M(\calA)=L\setminus\{0\}. 
$$
Choose a hyperplane $H$ containing $L$, then 
$\calA\cup\{H\}$ is free with exponents $(1, q-1, q)$. 
Again by using Theorem \ref{thm:addel}, we conclude that 
$\calA$ is free with exponents $(1, q-1, q-1)$. 
(2) can be proved similarly. 
\owari

We can summarize the results as follows. 

\begin{Cor}
(1) When $|\calA|\geq 2q$, $\calA$ is free if and only if $\chi(\calA, q)=0$, 
or equivalently, $M(\calA)=\emptyset$.  \\
(2) When $|\calA|=2q-1$, $\calA$ is free if and only if either 
$\chi(\calA, q)=0$ or 
$\chi(\calA, t)=(t-1)(t-q+1)^2$. \\
(3) When $|\calA|=2q-2$, $\calA$ is free if and only if either 
$\chi(\calA, q)=0$ or 
$\chi(\calA, q)=(t-1)(t-q+1)(t-q+2)$. 
\end{Cor}

\begin{Example}(\cite{zie-mat})\\
Let us fix a line $L\subset\bbF_3^3$, and define an arrangement $\calA_1$ by 
$$
\calA_1=\{H\in\calA_{all}(\bbF_3^3)\ |\ H\nsupseteq L\}.
$$
We can easily seen that $|\calA_1|=9>2\times 3$. Since 
$M(\calA_1)\neq\emptyset$, $\calA_1$ is not free. 
However, in \cite{zie-mat}, Ziegler proved the following. 
Let $\calA_2$ be an arrangement in $\bbK^3$ satisfying 
$L(\calA_1)\cong L(\calA_2)$. If the characteristic of the field $\bbK$ 
is not $3$, then $\calA_2$ is free with exponents $(1,4,4)$. 
\end{Example}

\medskip 

\noindent
{\bf Acknowledgements. }
This work was done 
while the special semester ``Hyperplane Arrangements and Application'' 
at the MSRI, Berkeley. 
The author thanks to the institute for their hospitality, and to the 
organizers for the wonderful semester. He also deeply thanks to 
Y. Kawahara, 
H. Terao, M. Wakefield and S. Yuzvinsky for many helpful discussions.


\begin{thebibliography}{99}

\bibitem[Ath]{ath-cha}
C. A. Athanasiadis,
Characteristic polynomials of subspace arrangements and finite fields,
Adv. Math. {\bf 122}(1996), 193--233.

\bibitem[CR]{cra-rot}
H. Crapo, G. -C. Rota,
On the Foundations of Combinatorial Theory:
Combinatorial Geometries, preliminary edition, MIT Press,
Cambridge, MA, 1970

\bibitem[OT]{orl-ter}
P. Orlik, H. Terao,
Arrangements of hyperplanes.
Grundlehren der Mathematischen Wissenschaften, 300.
Springer-Verlag, Berlin, 1992

\bibitem[Sa]{sai-log}
K. Saito,
Theory of logarithmic differential forms and logarithmic vector fields.
J. Fac. Sci. Univ. Tokyo, Sect. IA, Math. {\bf 27}(1980), 265--291.

\bibitem[Te1]{ter-arr}
H. Terao,
Arrangements of hyperplanes and their freeness. I.
J. Fac. Sci. Univ. Tokyo Sect. IA Math. {\bf 27} (1980), no. 2, 293--312.

\bibitem[Te2]{ter-free}
H. Terao,
Generalized exponents of a free arrangement of hyperplanes and
Shepherd-Todd-Brieskorn formula.
Invent. Math. {\bf 63} (1981), no. 1, 159--179.

\bibitem[Zi]{zie-mat}
G. Ziegler,
Matroid representations and free arrangements.
Trans. Amer. Math. Soc. {\bf 320}(1990), 525--541.


\end{thebibliography}
\end{document}